\documentclass[12pt]{article}
\usepackage{amssymb}
\usepackage{amsthm}
\usepackage{graphicx}
\usepackage{graphics}
\usepackage{amsthm}

\input pictex.tex
\newcommand{\grad}{\mathrm{o}}

\theoremstyle{definition}

\input epsf

\begin{document}

\date{}
\author{Andr\'es Navas}

\title{Pithagoras' theorem via equilateral triangles}
\maketitle

\vspace{-0.5cm}

\thispagestyle{empty}
\pagestyle{empty}

Pythagoras' theorem is perhaps the most important theorem in basic Mathematics, and many proofs are known. 
Besides Pythagoras', there is one in the Jiuzhang Suanshu, an ancient Chinese text which is somewhat contemporary to Pythagoras; see [1]. 
Another proof is attributed to Leonardo da Vinci [2], and there is also a proof by the USA president James Gardfield [3]. The reader may 
consult [4] for a list of 116 different proofs of the theorem, with discussion and links to many other ones.

Here we propose still another proof which is nonstandard as we do not use neither squares nor similarity of triangles. 
To motivate it, let $a,b$ be the cathetus of a right-angled triangle and $c$ its hypothenuse. Then from the equality 
$$a^2 + b^2 = c^2$$
it readily follows that the area of a regular $n$-gon of side $c$ is the sum of the areas of two regular $n$-gons of sides $a$ and $b$, respectively. 
To see this, it suffices to multiply both sides of this equality by an appropriate constant (namely, the area of the regular $n$-gon of side length 1). 
Actually, this statement for any given $n \geq 3$ also implies (hence it is equivalent to) Pythagoras' theorem just by reversing this argument.

\begin{figure}[ht!]
\centering
\includegraphics[width=36mm]{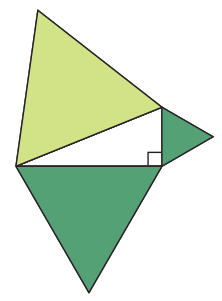}
\qquad
\includegraphics[width=40mm]{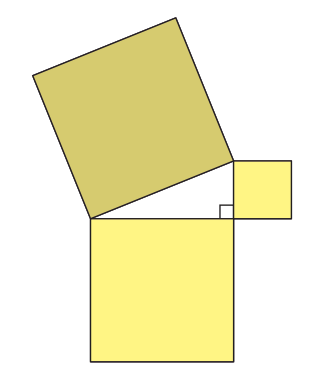}
\quad
\includegraphics[width=40mm]{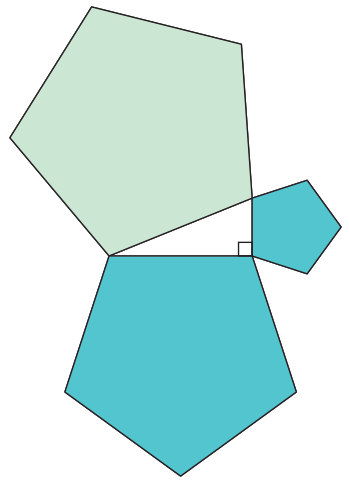}
%\caption{On the left, Nepal's flag; on the right, Togo's flag (for which the width-height ratio equals the golden mean). \label{nepal-togo}}
\end{figure}

In [5], Kassie Smith asks for a more geometric proof of all of this. Of course, in this framework, Wallace-Bolyai-Gerwien's decomposition 
theorem applies [6], hence there should be a decomposition of the two small $n$-gons into polygonal pieces that, after rearrangement, 
yield the larger one. Nevertheless, the number of pieces may be very large.

Below, we provide a short geometric proof of the statement above for equilateral triangles by using a beautiful and seemingly unnoticed 
configuration. To start with, let $A,B,C$ be the vertices of the triangle. Rotate the triangle $ABC$ counter-clockwise in $60^\grad$ at 
the vertex $A$, and clockwise also in $60^\grad$ at $B$. Call $C_1,B_1$ and $C_2,A_2$ the image of the vertices, as shown in 
the picture below. Notice that $BA_2 = c = AB_1$ and $\angle ABA_2 = \angle BAB_1 = 60^\grad$. Hence, $A_2$ and $B_1$ 
coincide, and denoting this point by $D$, the vertices $A,B$ and $D$ determine an equilateral triangle of side length $c$. 

\newpage

\begin{figure}[ht!]
\centering
\includegraphics[width=65mm]{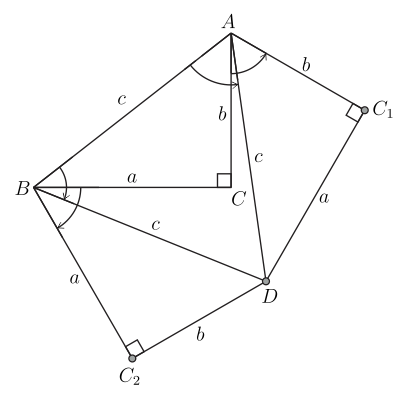}
%\caption{On the left, Nepal's flag; on the right, Togo's flag (for which the width-height ratio equals the golden mean). \label{nepal-togo}}
\end{figure}

Also, notice that $BCC_2$ and $ACC_1$ are equilateral triangles of side length $a$ and $b$, respectively. Moreover, 
triangles $BC_2D$ and $AC_1D$ are both congruent to triangle $BCA$. 

\begin{figure}[ht!]
\centering
\includegraphics[width=65mm]{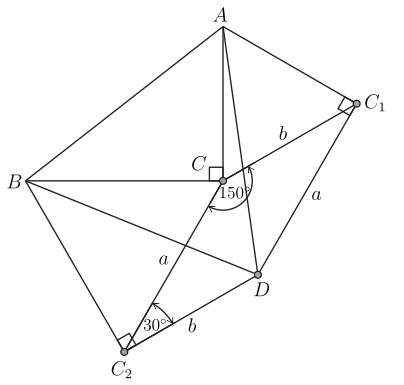}
%\caption{On the left, Nepal's flag; on the right, Togo's flag (for which the width-height ratio equals the golden mean). \label{nepal-togo}}
\end{figure}

Now, referring to the areas, we have 
$$ABC_2DC_1 = ABD + BC_2D + AC_1D = ACC_1 + BCC_2 + BCA + C_2 D C_1 C.$$
In pictures,

\begin{figure}[ht!]
\centering
\includegraphics[width=150mm]{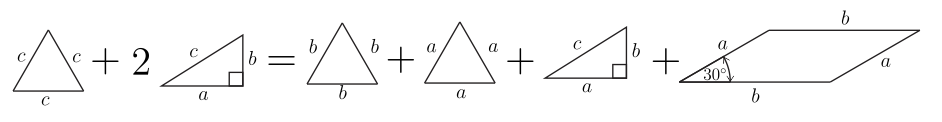}
%\caption{On the left, Nepal's flag; on the right, Togo's flag (for which the width-height ratio equals the golden mean). \label{nepal-togo}}
\end{figure}

\noindent which yields 

\begin{figure}[ht!]
\centering
\includegraphics[width=150mm]{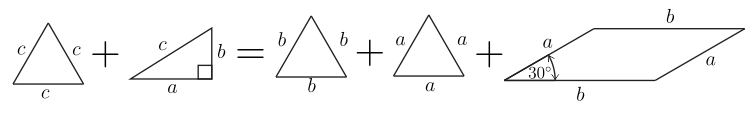}
%\caption{On the left, Nepal's flag; on the right, Togo's flag (for which the width-height ratio equals the golden mean). \label{nepal-togo}}
\end{figure}

\newpage

\noindent We are hence left to prove that 

\begin{figure}[ht!]
\centering
\includegraphics[width=90mm]{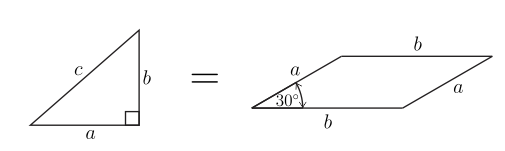}
%\caption{On the left, Nepal's flag; on the right, Togo's flag (for which the width-height ratio equals the golden mean). \label{nepal-togo}}
\end{figure}

\noindent To do this, notice that $C_2 D C_1 C$ is a parallelogram of sides $a$ and $b$. Besides, 
as $\angle BCA = 90^\grad$, we have that  $\angle C_1 C C_2$ must equal $120^\grad$, hence  
$\angle CC_2D = \angle DC_1C = 30^\grad$. Thus, 
$$C_1 C C_2 D = ab \sin (30^\grad) = \frac{1}{2} ab = BCA,$$
as desired.

\vspace{0.35cm}

\noindent{\bf Remark.} In the last step above, one can certainly avoid the use of trigonometry by just looking at the pictures below...

\vspace{-0.4cm}

\begin{figure}[ht!]
\centering
\includegraphics[width=120mm]{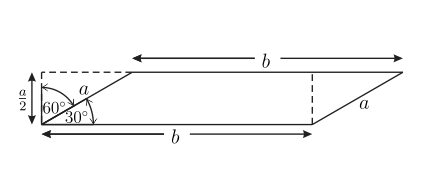}
%\caption{On the left, Nepal's flag; on the right, Togo's flag (for which the width-height ratio equals the golden mean). \label{nepal-togo}}
\end{figure}

\vspace{-0.4cm}

\begin{figure}[ht!]
\centering
\includegraphics[width=120mm]{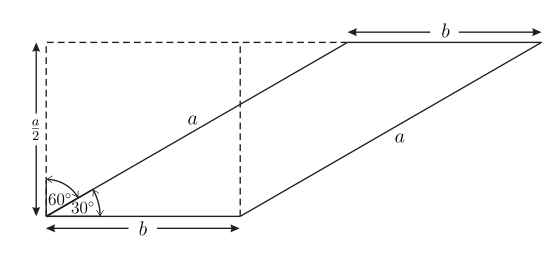}
%\caption{On the left, Nepal's flag; on the right, Togo's flag (for which the width-height ratio equals the golden mean). \label{nepal-togo}}
\end{figure}

%%%%%%%%%%%%%%%%%%%%%%%%%%%%%%%%%%%%%%%%%%%%%%%%%%%%%%%%%%%%%%%%%%%%%%%%%%%%%%%%%%%%%%%%%

\vspace{-0.5cm}

\noindent{\small  REFERENCES}

\vspace{0.2cm}

\begin{footnotesize}

%\begin{thebibliography}{Dillo 83}

%\bibitem{1} 
\noindent [1] https:$/\!/$math.temple.edu/$\sim$zit/Zitarelli/Pythag$_{-}$Chinese.pdf. 

%\bibitem{2} 
\noindent [2] http:$/\!/$www.rzuser.uni-heidelberg.de/$\sim$hb3/publ/vinci-0.pdf. 

%\bibitem{3} 
\noindent [3] http:$/\!/$math.kennesaw.edu/$\sim$sellerme/sfehtml/classes/math1112/garfieldpro.pdf. 

%\bibitem{4} 
\noindent [4] http:$/\!/$www.cut-the-knot.org/pythagoras/.

%\bibitem{5} 
\noindent [5] http:$/\!/$jwilson.coe.uga.edu/EMAT6680Fa2012/Smith/6690/pythagorean$\%$20theorem/KLS$_{-}$Pythagorean$_{-}$Theorem.html.

\noindent [6] https:$/\!/$en.wikipedia.org/wiki/Wallace$\%$E2$\%$80$\%$93Bolyai$\%$E2$\%$80$\%$93Gerwien$_{-}$theorem.

%\end{thebibliography}

\vspace{0.35cm}

\noindent Andr\'es Navas (andres.navas@usach.cl)\\ 

\noindent Dpto. de Matem\'atica y Ciencia de la Computaci\'on, Universidad de Santiago de Chile\\ 

\noindent Alameda 3363, Santiago, Chile

\end{footnotesize}

\end{document}